\documentclass{ifacconf}

\usepackage{graphicx}      
\usepackage{natbib}        

\usepackage{amssymb}
\usepackage{amsmath,amsfonts}
\newtheorem{definition}{Definition}
\newtheorem{property}{Property}
\newtheorem{lemma}{Lemma}
\newtheorem{theorem}{Theorem}
\newtheorem{remark}{Remark}
\newtheorem{assumption}{Assumption}
\usepackage{graphicx}
\usepackage{float}
\usepackage{subfigure}
\usepackage{adjustbox}
\begin{document}
\begin{frontmatter}

\title{Interval Estimation for Bounded Jacobian Nonlinear Systems by Zonotope Analysis\thanksref{footnoteinfo}} 

\thanks[footnoteinfo]{This work was supported by the National Natural Science Foundation of China under Grant No. 62373125.}

\author[First]{Chi Xu} 
\author[First]{Zhenhua Wang} 
\author[fouth]{ Nacim Meslem}
\author[fifth]{Tarek Raïssi}
\author[Third]{Yacine Chitour} 
\address[First]{ Department of Control Science and Engineering, Harbin Institute of Technology, Harbin, 150001 P. R. China, (e-mail:  chi.xu@stu.hit.edu.cn, zhenhua.wang@hit.edu.cn).}
\address[fouth]{ University Grenoble Alpes, CNRS, Grenoble INP, GIPSA-lab, F-38000 Grenoble, France , (e-mail:  nacim.meslem@grenoble-inp.fr)}
\address[fifth]{Conservatoire National des Arts et Métiers (CNAM), 292 Rue Saint-Martin, 75141 Paris, France,(e-mail: author@snu.ac.kr)}
\address[Third]{ Université Paris-Saclay, CNRS, CentraleSupélec, Laboratoire des Signaux et Systèmes, 91190 Gif sur-Yvette, France, and also with Fédération de Mathématiques de CentraleSupélec, 91190 Gif-sur-Yvette, France (e-mail: yacine.chitour@centralesupelec.fr)}
\begin{abstract}                
This paper introduces an interval state estimation method for discrete-time bounded Jacobian nonlinear systems allying Luenberger-like observer with zonotope set computation.
First, a robust observer is designed to obtain bounded-error and point-estimation with a peak-to-peak performance.
This allows one to cope with the unknown but bounded process disturbance and measurement noise. 
Then, based on the stable dynamics of the observation error, tight interval estimation is obtained by applying zonotope set computation and analysis.
To sum up, a comprehensive interval estimation algorithm is proposed by integrating the robust (in the sense of peak-to-peak performance index) point-valued estimate with the feasible zonotope set of the estimation error.
Numerical simulation tests are conducted to assess the effectiveness of the proposed approach.
\end{abstract}


\end{frontmatter}

\section{INTRODUCTION}
The knowledge of the internal state of dynamics systems is essential in many engineering applications, including systems monitoring, fault detection, and state feedback control \cite{wang2023model}.
Due to uncertainties such as process disturbances and measurement noise, the estimated state of a system may not align perfectly with its actual one. 
A solution to solve this problem (to get a reliable estimation) can be obtained by assuming that the uncertain parts of the system are unknown but bounded and then using interval computation techniques to design set-valued state estimators. 
It is worth pointing out that compared to the stochastic assumption, that is the existence of probability distribution functions describing the system uncertainties, the unknown but bounded assumption is more natural and easier to meet in real world systems.

In the literature, interval estimation techniques can be divided into two main categories, interval observer-based methods \cite{li2019interval,tahir2021synthesis,khajenejad2022h} and set-based methods \cite{wang2018zonotopic,zhang2021ellipsoid,wang2024interval}.  
Using the monotone systems theory, the interval observer-based methods aim to design two sub-observers to ensure the stability and positivity of the estimation error, and provide two bounds framing the system state \cite{raissi2011interval}.
However, sometimes the design condition of interval observers is too restrictive to find suitable gain matrices in the sense of estimation performance, which limits their application. 
Although using coordinate transformation can relax the design condition, the transformation  introduces additional conservatism and reduces the estimation accuracy. 

Unlike interval observers, the main idea of set-based methods is to construct, step-by-step, geometrical sets that contain the reachable state of the system \cite{althoff2021comparison}. Ellipsoids, polytopes and zonotopes are commonly used as geometrical sets. 
Owing of their geometrical and computational properties, zonotope-based methods offer a good trade-off between estimation accuracy and computational efficiency \cite{tang2019interval}.
For nonlinear discrete-time systems, the authors of \cite{alamo2005guaranteed} have proposed a zonotope-based interval state estimator based on an adaptation of the prediction/correction principle to the set-membership paradigm. Unluckily, the convergence rate of the estimation error of this approach is not tunable and there is no proof of its stability  \cite{efimov2013interval}. For the case of linear time-invariant systems (LTI), the authors of \cite{tang2019interval} have developed an efficient two-step set-valued estimation method by coupling a point-wise observer scheme with zonotope analysis of the estimation error. Unfortunately, the application scope of this method is limited to the case of linear systems and cannot be directly extended to the case of nonlinear systems.

To overcome the limitations of the aforementioned works, we propose an interval estimation method that combines a point estimation observer with zonotope analysis of the estimation error for bounded Jacobian nonlinear systems. The main contributions of this work are summarized as follows:
\begin{enumerate}
	\item A bounded-error observer is designed to provide a point estimation with a peak-to-peak performance. In fact, the proposed design method for the observer gain allows one attenuating the negative effects of the unknown but bounded state disturbance and measurement noise on the estimation error.
	\item By using the available bounds on the state disturbances and measurement noise, online zonotope computations are applied to design a comprehensive set-valued state estimator. Compared with the existing methods [3], the proposed state estimator provides results that are more accurate.
\end{enumerate}

The remainder of this paper is organized as follows. 
Section \ref{sec2} introduces some preliminaries and problem formulation. 
In Section \ref{sec3}, a comprehensive interval estimation method is introduced based on a bounded-error observer design and a zonotope analysis. To illustrate the effectiveness of the proposed method, numerical simulation results are given and commented in Section \ref{sec4}.
A conclusion and some perspectives of this work are presented in Section \ref{sec5}.

\section{PRELIMINARIES AND PROBLEM FORMULATION}
\label{sec2}
\subsection{Preliminaries}
\textbf{Notation}: $\mathbb{R}^n$ denotes the $n$-dimensional Euclidean space. $\mathbb{R}^{m\times n}$ is the set of matrices of size $m$ by $n$. 
$\oplus$ stands for the Minkowski sum. 
$(\cdot)^{\mathrm{T}}$ represents the transpose of a matrix/vector.
In a symmetric matrix, asterisk $*$ stands for the symmetric term.
$\min(\cdot)$, $\max(\cdot)$, $\leq,<,\geq,>$ and $|\cdot|$ should be understood elementwise. 
For a vector $\boldsymbol{z} \in \mathbb{R}^n$, $\boldsymbol{z}(i)$ is its $i$-th element, $\|\boldsymbol{z}\|_1 = \sum_{i}^{n}|\boldsymbol{z}(i)|$, $\|\boldsymbol{z}\|_2 = \sqrt{\boldsymbol{z}^\mathrm{T}\boldsymbol{z}}$,  $\mathrm{diag}(\boldsymbol{z})$ represents a diagonal matrix with the entries of $\boldsymbol{z}$ as its diagonal elements.
An interval  $[\boldsymbol{z}]\triangleq[\underline{\boldsymbol{z}},\overline{\boldsymbol{z}}]$ is a set of vectors $\boldsymbol{z}$ satisfying $\underline{\boldsymbol{z}}\leq \boldsymbol{z}\leq \overline{\boldsymbol{z}}$.
For a matrix $\boldsymbol{M}\in\mathbb{R}^{n\times m}$, $\boldsymbol{M}(i,j)$ is its $(i,j)$-th element, $\boldsymbol{M}\succ \boldsymbol{0} (\prec \boldsymbol{0}) $ means that $\boldsymbol{M}$ is a positive (negative) definite matrix and an interval matrix $[\boldsymbol{M}]\triangleq[\underline{\boldsymbol{M}},\overline{\boldsymbol{M}}]$ is a set of matrices $\boldsymbol{M}$ which satisfy $\underline{\boldsymbol{M}}(i,j)\leq
\boldsymbol{M}(i,j)\leq \overline{\boldsymbol{M}}(i,j), i=1,\dots,n, j=1,\dots , m$. 
$\downarrow_{q}(\boldsymbol{M})$ is a zonotope reduction operator introduced in Appendix.
$\boldsymbol{I}_{n}$ denotes an ${n \times n}$ identity matrix and $\boldsymbol{0}$ is a zero matrix/vector with suitable dimensions.

\begin{definition}\label{zonotope}
	Given a vector $\boldsymbol{p}\in \mathbb{R}^{n}$ and a generator matrix $\boldsymbol{M}\in \mathbb{R}^{n \times m}$,
	an $m$-order zonotope $\mathcal{X} \subseteq \mathbb{R}^{n} $ is defined as:
	\begin{equation}
		\begin{adjustbox}{max width=0.9\linewidth}
			$
			\mathcal{X} = \boldsymbol{p} \oplus \boldsymbol{M}\mathcal{B}^{m}
			=\left\{ \boldsymbol{p}  + \boldsymbol{M}\xi, \xi \in \mathcal{B}^{m}, \mathcal{B}^{m} = [-1,1]^{m}\right\}
			$
		\end{adjustbox}
	\end{equation}
	For simplicity, in what follows, a zonotope is represented as $ \mathcal{X} = \langle \boldsymbol{p}, \boldsymbol{M} \rangle $.
\end{definition}

\begin{property}\cite{tang2020set}\label{Pro1}
	Given a matrix $\boldsymbol{L}\in\mathbb{R}^{l\times n}$, a vector $\boldsymbol{z}\in \mathbb{R}^{n}$ and zonotopes $\langle \boldsymbol{p}, \boldsymbol{M} \rangle,\langle \boldsymbol{p}_1, \boldsymbol{M}_1 \rangle,\langle \boldsymbol{p}_2, \boldsymbol{M}_2 \rangle \subseteq \mathbb{R}^{n}$ one has:
	\begin{align}
		\boldsymbol{L}\langle \boldsymbol{p}, \boldsymbol{M} \rangle\oplus\boldsymbol{z} &= \langle \boldsymbol{Lp}+\boldsymbol{z}, \boldsymbol{LM} \rangle \\
		\boldsymbol{x} \in \langle \boldsymbol{p}, \boldsymbol{M} \rangle &\subseteq  [\underline{\boldsymbol{x}},\overline{\boldsymbol{x}}]\\
		\langle \boldsymbol{p}, \boldsymbol{M} 	\rangle &\subseteq \langle
		\boldsymbol{p}, \downarrow_{q}(\boldsymbol{M}) \rangle \\
		\langle \boldsymbol{p}_1, \boldsymbol{M}_1 \rangle \oplus \langle \boldsymbol{p}_2, \boldsymbol{M}_2 \rangle &= \langle \boldsymbol{p}_1 +\boldsymbol{p}_2 , [\boldsymbol{M}_1 \!\quad\! \boldsymbol{M}_2] \rangle 
	\end{align}
	where $\boldsymbol{p},\boldsymbol{p}_1,\boldsymbol{p}_2 \in \mathbb{R}^n$, $\boldsymbol{M}\in \mathbb{R}^{n \times m},\boldsymbol{M}_1\in \mathbb{R}^{n \times m_1},\boldsymbol{M}_2\in \mathbb{R}^{n \times m_2}$.
	$\downarrow_{q}(\cdot)$ is a zonotope reduction operator and $q$ is a chosen integer satisfying $n<q<m$. 
	Moreover, the bounds of $\boldsymbol{x}$ can be calculated by:
	\begin{equation}\label{eq6}
		\left\lbrace \begin{array}{c}
			\underline{\boldsymbol{x}}(i) = \boldsymbol{p}(i) - \|\boldsymbol{M}(i,:)\|_1 \\
			\overline{\boldsymbol{x}}(i) = \boldsymbol{p}(i) + \|\boldsymbol{M}(i,:)\|_1
		\end{array}\right.,\quad i = 1,\dots,n  
	\end{equation}
\end{property}

\begin{lemma}\cite{de1999new}\label{lem1}
	Given a symmetric matrix $\boldsymbol{M}\succ \boldsymbol{0}$ and matrices $\boldsymbol{N,Q,U}$, the following conditions are equivalent:
	\begin{align}
		&\boldsymbol{Q}^{\mathrm{T}}\boldsymbol{MQ} - \boldsymbol{N}\prec \boldsymbol{0}\\
		&\left\lbrack \begin{array}{cc}
			-\boldsymbol{N} & *\\
			\boldsymbol{UQ} &  \boldsymbol{M}-\boldsymbol{U}-\boldsymbol{U}^{\mathrm{T}} 
		\end{array}\right\rbrack \prec \boldsymbol{0}
	\end{align}
\end{lemma}

\begin{lemma}\cite{alamo2005guaranteed}\label{lemma2}
	Given a vector $\boldsymbol{p}\in \mathbb{R}^{n}$ and an interval matrix $[\boldsymbol{M}] = [\underline{\boldsymbol{M}},\overline{\boldsymbol{M}}] \subseteq \mathbb{R}^{n \times m} $, 
	a group of zonotopes is denoted as $\mathcal{X} = \boldsymbol{p} \oplus [\boldsymbol{M}]\mathcal{B}^m$. 
	A zonotope inclusion is a zonotope $\lozenge \mathcal{X}$ that contains $\mathcal{X}$ and computed as follows:
	\begin{equation}
		\begin{aligned}
			\mathcal{X} \subset \lozenge \mathcal{X} &= \boldsymbol{p} \oplus 
			\left[
			\boldsymbol{D^c_{[M]}}, \boldsymbol{D^s_{[M]}}
			\right]\mathcal{B}^{m+n} \\
			& = \langle
			\boldsymbol{p}, 	\left[
			\boldsymbol{D^c_{[M]}}, \boldsymbol{D^s_{[M]}}
			\right] \rangle 
		\end{aligned}
	\end{equation}
	where $\boldsymbol{D}^c_{[\boldsymbol{M}]} \in \mathbb{R}^{n\times m}$, $\boldsymbol{D}^s_{[\boldsymbol{M}]}\in \mathbb{R}^{n\times n}$ are matrices related to $[\boldsymbol{M}]$ which are defined by:
	\begin{equation}\label{eq10}
		\begin{aligned}
			\boldsymbol{D}^c_{[\boldsymbol{M}]} &= \frac{\underline{\boldsymbol{M}}+\overline{\boldsymbol{M}}}{2}, \ \boldsymbol{R}^c_{[\boldsymbol{M}]} = \frac{\overline{\boldsymbol{M}}-\underline{\boldsymbol{M}}}{2},\\
			\boldsymbol{s}_{[\boldsymbol{M}]}(i) &= \| \boldsymbol{R}^c_{[\boldsymbol{M}]}(i,:)\|_1,\ i,\dots,n,\\
			\boldsymbol{D}^s_{[\boldsymbol{M}]} &= \mathrm{diag}(\boldsymbol{s}_{[\boldsymbol{M}]}).
		\end{aligned}
	\end{equation}
	
\end{lemma}

\subsection{Problem statement}
Consider a discrete-time nonlinear system described by 
\begin{align}\label{eq11}
	\left\lbrace 
	\begin{aligned}
		\boldsymbol{x}_{k+1}&=\boldsymbol{A}\boldsymbol{x}_{k}+\boldsymbol{Bu}_{k}+\boldsymbol{f}(\boldsymbol{x}_k)+\boldsymbol{D}_1\boldsymbol{w}_{k}\\
		\boldsymbol{y}_k &= \boldsymbol{Cx}_{k}+\boldsymbol{D}_2\boldsymbol{v}_{k}
	\end{aligned}
	\right. 
\end{align} 
where $\boldsymbol{x}_k \in \mathbb{R}^{n_x}$, $\boldsymbol{y}_k \in\mathbb{R}^{n_y}$, $\boldsymbol{u}_k\in \mathbb{R}^{n_u}$, $\boldsymbol{w}_k\in \mathbb{R}^{n_w}$ and $\boldsymbol{v}_k\in \mathbb{R}^{n_v}$ are the system state, the measurement output, the control input,  the process disturbance, and the measurement noise, respectively. $\boldsymbol{A,B,C},\boldsymbol{D}_1,\boldsymbol{D}_2$ are known matrices with suitable dimensions. $\boldsymbol{f}(\cdot): \mathbb{R}^{n_x} \rightarrow \mathbb{R}^{n_x} $ is a first-order  differentiable nonlinear function with respect to $\boldsymbol{x}$.

\begin{assumption}
	Let $\mathcal{X} \subseteq \mathbb{R}^{n_x}$ be the feasible domain of $\boldsymbol{x}$. 
	There exist known matrices $\underline{\boldsymbol{J}}, \overline{\boldsymbol{J}} \in \mathbb{R}^{n_x \times n_x}$ such that for all $\boldsymbol{x} \in \mathcal{X}$, the Jacobian of $\boldsymbol{f}(\cdot)$ satisfies
	\begin{align}
		\underline{\boldsymbol{J}} \leq \frac{\partial \boldsymbol{f}}{\partial \boldsymbol{x}}(\boldsymbol{x}) \leq \overline{\boldsymbol{J}}.
	\end{align}
\end{assumption}
Under Assumption 1, for any $\boldsymbol{a}, \boldsymbol{b} \in \mathcal{X}$, there exists a matrix $\boldsymbol{J} \in [\underline{\boldsymbol{J}}, \overline{\boldsymbol{J}}]$ such that
\begin{align}
	\label{eq13}\Delta \boldsymbol{f}^{\mathrm{T}}(\boldsymbol{a}-\boldsymbol{b}) &\leq (\boldsymbol{a}-\boldsymbol{b})^{\mathrm{T}}\boldsymbol{J}(\boldsymbol{a}-\boldsymbol{b}), \\
	\label{eq14}\Delta \boldsymbol{f}^{\mathrm{T}}\Delta \boldsymbol{f} &\leq (\boldsymbol{a}-\boldsymbol{b})^{\mathrm{T}}\boldsymbol{J}^{\mathrm{T}}\boldsymbol{J}(\boldsymbol{a}-\boldsymbol{b}),
\end{align}
where $\Delta \boldsymbol{f} = \boldsymbol{f}(\boldsymbol{a}) - \boldsymbol{f}(\boldsymbol{b})$. If $\boldsymbol{J}^\mathrm{T}\boldsymbol{J}=\gamma_f^2\boldsymbol{I}$, inequality (\ref{eq14}) aligns with the standard Lipschitz condition. If $\boldsymbol{J}=\rho_f \boldsymbol{I}$, the (\ref{eq13}) aligns with the one-sided Lipschitz condition. Thus, Assumption 1 holds for a large class of systems.

\begin{assumption}
	The process disturbance $\boldsymbol{w}_k$, measurement noise $\boldsymbol{v}_k$ and initial state $\boldsymbol{x}_0$ satisfy unknown but bounded conditions as follows:
	\begin{align}\label{eq15}
		\forall \ k\geq 0,\
		|\boldsymbol{w}_k| \leq \overline{\boldsymbol{w}},\ |\boldsymbol{v}_k| \leq \overline{\boldsymbol{v}}, \ |\boldsymbol{x}_0-\boldsymbol{p}_0|\leq \overline{\boldsymbol{x}}_0,
	\end{align}
	where $\overline{\boldsymbol{w}}\in \mathbb{R}^{n_w},\overline{\boldsymbol{v}}\in \mathbb{R}^{n_v}$, $\overline{\boldsymbol{x}}_0 \in \mathbb{R}^{n_x}$ and $\boldsymbol{p}_0 \in \mathbb{R}^{n_x}$ are known. 
	
	Based on (\ref{zonotope}),  (\ref{eq15}) can be reformulated as 
	\begin{align}
		\begin{aligned}
			\boldsymbol{w}_k &\in \mathcal{W} = \langle \boldsymbol{0}, \boldsymbol{M}_w \rangle, 	
			\boldsymbol{v}_k \in \mathcal{V} = \langle \boldsymbol{0}, \boldsymbol{M}_v \rangle, \\
			\boldsymbol{x}_0 &\in \mathcal{X}_0 = \langle \boldsymbol{p}_0, \boldsymbol{M}_{\boldsymbol{x},0} \rangle,
		\end{aligned}
	\end{align}
	where $\boldsymbol{M}_w=\mathrm{diag}(\overline{\boldsymbol{w}})$, $\boldsymbol{M}_v=\mathrm{diag}(\overline{\boldsymbol{v}})$, $\boldsymbol{M}_{\boldsymbol{x},0}=\mathrm{diag}(\overline{\boldsymbol{x}}_0)$.
\end{assumption}

This paper aims at designing a set-valued state estimator for system (11) that computes, in real-time, a tight interval $[\underline{\boldsymbol{x}},\overline{\boldsymbol{x}}]$ that contains the actual state vector of the system.

\section{MAIN RESULT}\label{sec3}
In this section, a robust observer is designed to provide point-estimation in the bounded-error context with peak-to-peak performance index. Then, a comprehensive interval estimation method is proposed by combining the  point-estimate and zonotope analysis of the reachable set of the estimation error.
\subsection{An offline Bounded-error Observer Design}
Consider the following observer
\begin{align}\label{eq17}
	\hat{\boldsymbol{x}}_{k+1} = \boldsymbol{A}\hat{\boldsymbol{x}}_k + \boldsymbol{Bu}_k + \boldsymbol{f}(\hat{\boldsymbol{x}}_k)+\boldsymbol{L}(\boldsymbol{y}_k-\boldsymbol{C}\hat{\boldsymbol{x}}_k)
\end{align}
where $\hat{\boldsymbol{x}} \in \mathbb{R}^{n_x}$ and $ \boldsymbol{L} \in \mathbb{R}^{n_x\times n_y}$ are the estimated state and observer gain.

The estimation error is defined as: 
\begin{align}\label{eq18}
	\boldsymbol{e}_k = \boldsymbol{x}_k - \hat{\boldsymbol{x}}_k.
\end{align}
By substituting (\ref{eq11}) and (\ref{eq17}) into (\ref{eq18}), we have
\begin{align}\label{eq19}
	\boldsymbol{e}_{k+1}= \tilde{\boldsymbol{A}}\boldsymbol{e}_{k} + \Delta \boldsymbol{f}_k + \boldsymbol{D}_1\boldsymbol{w}_k-\boldsymbol{L}\boldsymbol{D}_2\boldsymbol{v}_k
\end{align}
with $\tilde{\boldsymbol{A}} = \boldsymbol{A}-\boldsymbol{LC}$, $\Delta \boldsymbol{f}_k = \boldsymbol{f}({\boldsymbol{x}}_k) - \boldsymbol{f}(\hat{\boldsymbol{x}}_k)$. 

\begin{theorem}
	For a given constant \( 0 < \lambda < 1 \), suppose there exist constants \( a > 0 \), \( b > 0 \), \( \gamma > 0 \), \( \mu_w > 0 \), \( \mu_v > 0 \), 
	a positive definite matrix \(\boldsymbol{P} \succ \boldsymbol{0} \in \mathbb{R}^{n_x \times n_x}\), 
	matrices \(\boldsymbol{W} \in \mathbb{R}^{n_x \times n_y}\), \(\boldsymbol{U} \in \mathbb{R}^{n_x \times n_y}\), 
	such that for all \(\boldsymbol{J} \in [\underline{\boldsymbol{J}}, \overline{\boldsymbol{J}}]\), the following inequalities hold:
	\begin{align}
		\label{eq20}    
		\begin{bmatrix}
			\boldsymbol{\Phi}_{11} & * & * & * & * & *\\
			-a\boldsymbol{I} & -b\boldsymbol{I}_{n_x} &* &* & * & *\\
			\boldsymbol{0} & \boldsymbol{0} & -\gamma\boldsymbol{I}_{n_w} & * & * & *\\
			\boldsymbol{0} & \boldsymbol{0} & \boldsymbol{0} & -\gamma\boldsymbol{I}_{n_v} & * & *\\
			\boldsymbol{\Phi}_{51} & \boldsymbol{P} &\boldsymbol{PD}_1 & -\boldsymbol{WD}_2 & -\boldsymbol{P} & *\\
			\boldsymbol{UJ} & \boldsymbol{0} & \boldsymbol{0} & \boldsymbol{0}& \boldsymbol{0} & \boldsymbol{\Phi}_{66}
		\end{bmatrix} \prec \boldsymbol{0}, \\
		\label{eq21}    
		\begin{bmatrix}
			\lambda\boldsymbol{P} & *& * &* \\
			\boldsymbol{0}& (\mu_w - \gamma)\boldsymbol{I}_{n_w} & *& *\\
			\boldsymbol{0}& \boldsymbol{0}& (\mu_v - \gamma)\boldsymbol{I}_{n_v} & *\\
			\boldsymbol{I}_{n_x} & \boldsymbol{0}& \boldsymbol{0}& \mu \boldsymbol{I}_{n_x}
		\end{bmatrix} \succ \boldsymbol{0},
	\end{align}
	where 
	
	$\boldsymbol{\Phi}_{11} = \lambda\boldsymbol{P} - \boldsymbol{P} + a(\boldsymbol{J} + \boldsymbol{J}^\mathrm{T}), \
	\boldsymbol{\Phi}_{51} = \boldsymbol{PA} - \boldsymbol{WC}, $
	
	$ \boldsymbol{\Phi}_{66} = b\boldsymbol{I} - \boldsymbol{U} - \boldsymbol{U}^\mathrm{T}, \
	\mu = \mu_w + \mu_v.$
	
	Then, the estimation error dynamics (\ref{eq19}) is input-to-state stable 
	and its solution satisfies the peak-to-peak performance bound:
	\begin{align}\label{eq23}
		\|\boldsymbol{e}_k \|_2 \leq \sqrt{\mu \left( \lambda(1-\lambda)^k V_0 + \mu_w \overline{\boldsymbol{w}}^2 + \mu_v \overline{\boldsymbol{v}}^2 \right)},
	\end{align}
	where \( V_0 = \boldsymbol{e}_0^\mathrm{T} \boldsymbol{P} \boldsymbol{e}_0 \) and \(\mu\) quantifies the attenuation level.
	
	Thus, optimal estimation performance can be obtained by solving the following minimization problem:
	\begin{align}\label{eq24}
		\min_{a, b, \gamma, \mu_w, \mu_v, \boldsymbol{P}, \boldsymbol{W}, \boldsymbol{U}} \mu \quad \text{s. t.} \quad (\ref{eq20}), (\ref{eq21}).
	\end{align}
	
	If problem (\ref{eq24}) is feasible, the observer gain matrix \(\boldsymbol{L}\) is given by
	\begin{align}
		\boldsymbol{L} = \boldsymbol{P}^{-1} \boldsymbol{W}.
	\end{align}
\end{theorem}

Proof:
	Firstly, using the positive definite matrix $\boldsymbol{P}\succ \boldsymbol{0}$, we define the following Lyapunov function
	\begin{align}
		V_k = \boldsymbol{e}_k^{\mathrm{T}}\boldsymbol{P}\boldsymbol{e}_k
	\end{align} 
	and its time variation is described by
	\begin{align}
		\Delta V_k = V_{k+1} - V_k
		= \boldsymbol{\zeta}_k^\mathrm{T}\boldsymbol{\Psi}\boldsymbol{\zeta}_k
	\end{align} 
	where 
	\begin{align}\nonumber
		\begin{aligned}
			\boldsymbol{\zeta}_k&=\left[ \begin{matrix}
				\boldsymbol{e}_k\\ \Delta\boldsymbol{f}_k \\ \boldsymbol{w}_k \\ \boldsymbol{v}_k
			\end{matrix}\right],\ 
			\boldsymbol{\Psi}=\left[\begin{matrix}
				\Psi_{11} & *  & *  & * \\
				\Psi_{21} & \Psi_{22}  & *  & * \\
				\Psi_{31} & \Psi_{32}  & \Psi_{33}  & * \\
				\Psi_{41} & \Psi_{42}  & \Psi_{43}  & \Psi_{44}
			\end{matrix} \right], \\ 
			\Psi_{11} &= \tilde{\boldsymbol{A}}^{\mathrm{T}}\boldsymbol{P}\tilde{\boldsymbol{A}}-\boldsymbol{P},\
			\Psi_{21} = \boldsymbol{P}\tilde{\boldsymbol{A}},\ \Psi_{22} = \boldsymbol{P},\\
			\Psi_{31} &= \boldsymbol{D}_1^{\mathrm{T}}\boldsymbol{P}\tilde{\boldsymbol{A}}, \ \Psi_{32} = \boldsymbol{D}_1^{\mathrm{T}}\boldsymbol{P},\
			\Psi_{33} =\boldsymbol{D}_1^{\mathrm{T}}\boldsymbol{P}\boldsymbol{D}_1,\\
			\Psi_{41} &=(-\boldsymbol{LD}_2)^{\mathrm{T}}\boldsymbol{P}\tilde{\boldsymbol{A}},\ 
			\Psi_{42} =(-\boldsymbol{LD}_2)^{\mathrm{T}}\boldsymbol{P},\\
			\Psi_{43} &=(-\boldsymbol{LD}_2)^{\mathrm{T}}\boldsymbol{PD}_1,\ \Psi_{44} =(-\boldsymbol{LD}_2)^{\mathrm{T}}\boldsymbol{P}(-\boldsymbol{LD}_2).\\
		\end{aligned}
	\end{align}
	
	Then,  using Lemma \ref{lem1} and substituting $\boldsymbol{W}=\boldsymbol{PL}$ into (\ref{eq20}), it yields 
	\begin{align} \label{eq28}
		\left[\begin{matrix}
			\boldsymbol{\Gamma}_{11} & * & * & * & * \\
			-a\boldsymbol{I} & -b\boldsymbol{I}_{n_x} &* &* & * \\
			\boldsymbol{0} & \boldsymbol{0} & -\gamma_l\boldsymbol{I}_{n_w} & * & * \\
			\boldsymbol{0} & \boldsymbol{0} & \boldsymbol{0} & -\gamma_l\boldsymbol{I}_{n_v} & * \\
			\boldsymbol{\Gamma}_{51}& \boldsymbol{P} &\boldsymbol{PD}_1 & -\boldsymbol{WD}_2 & -\boldsymbol{P} 
		\end{matrix} \right] \prec \boldsymbol{0}
	\end{align}
	where
	\begin{align}\nonumber
		\begin{aligned}
			\boldsymbol{\Gamma}_{11} &= \lambda_l\boldsymbol{P}-\boldsymbol{P}+a(\boldsymbol{J}+\boldsymbol{J}^\mathrm{T})+b\boldsymbol{J}^\mathrm{T}\boldsymbol{J},\\
			\boldsymbol{\Gamma}_{51} &=\boldsymbol{P}(\boldsymbol{A}-\boldsymbol{LC}).
		\end{aligned}
	\end{align}
	Using the Schur complement lemma \cite{boyd2004convex}, (\ref{eq28}) can be rewritten as
	\begin{align}\label{eq29}
		\boldsymbol{\Psi} + \left[ \begin{matrix}
			\Xi_{11} &*  &* &* \\
			-a\boldsymbol{I}& -b\boldsymbol{I}_{n_x} &* &* \\
			\boldsymbol{0} & \boldsymbol{0} & -\gamma_l\boldsymbol{I}_{nw} & *\\
			\boldsymbol{0} & \boldsymbol{0} & \boldsymbol{0} & -\gamma_l\boldsymbol{I}_{nv} \\
		\end{matrix}\right] \prec \boldsymbol{0}
	\end{align} 
	where
	$
	\Xi_{11}= \lambda_l\boldsymbol{P}+a(\boldsymbol{J}+\boldsymbol{J}^\mathrm{T})+b\boldsymbol{J}^\mathrm{T}\boldsymbol{J}
	$
	
	Pre- and post-multiplying (\ref{eq29}) with $\boldsymbol{\zeta}_k^\mathrm{T}$ and its transpose yields 
	\begin{align} \label{eq30}
		\Delta V_k + \lambda_l V_k - \gamma_l \boldsymbol{w}_k^\mathrm{T}\boldsymbol{w}_k-\gamma_l \boldsymbol{v}_k^\mathrm{T}\boldsymbol{v}_k + \alpha + \beta \leq 0
	\end{align}
	where
	\begin{align}\nonumber
		\begin{aligned}
			\alpha &=a\boldsymbol{e}^\mathrm{T}_k(\boldsymbol{J}+\boldsymbol{J}^\mathrm{T})\boldsymbol{e}_k
			-a\boldsymbol{\Delta f}_k^{\mathrm{T}}{\boldsymbol{e}}_k 
			-a {\boldsymbol{e}}_k^{\mathrm{T}}\boldsymbol{\Delta f}_k \\
			\beta  &=b {\boldsymbol{e}}_{k}^\mathrm{T}\boldsymbol{J}^\mathrm{T}\boldsymbol{J}{\boldsymbol{e}}_k - 	b \boldsymbol{\Delta f}_k^{\mathrm{T}}\boldsymbol{\Delta f}_k
		\end{aligned}
	\end{align}
	At this stage one can claim that if (\ref{eq20}) holds, (\ref{eq30}) is satisfied.
	
	Subsequently, based on (\ref{eq13}) and $a>0$, it follows that 
	\begin{align}
		\Delta \boldsymbol{f}^{\mathrm{T}}_k(\boldsymbol{x}_k-\hat{\boldsymbol{x}}_k)&\leq (\boldsymbol{x}_k-\hat{\boldsymbol{x}}_k)^{\mathrm{T}}\boldsymbol{J}(\boldsymbol{x}_k-\hat{\boldsymbol{x}}_k) \nonumber\\
		\Delta \boldsymbol{f}^{\mathrm{T}}_k\boldsymbol{e}_k&\leq\boldsymbol{e}_k^\mathrm{T}\boldsymbol{J}\boldsymbol{e}_k \nonumber\\
		\boldsymbol{\Delta f}_k^{\mathrm{T}}{\boldsymbol{e}}_k 
		+{\boldsymbol{e}}_k^{\mathrm{T}}\boldsymbol{\Delta f}_k &\leq \boldsymbol{e}^\mathrm{T}_k(\boldsymbol{J}+\boldsymbol{J}^\mathrm{T})\boldsymbol{e}_k \nonumber \\
		\label{eq31}\alpha &\geq 0
	\end{align}
	Similarly, based on (\ref{eq14}) and $b>0$, we have 
	\begin{align}
		\Delta \boldsymbol{f}^{\mathrm{T}}_k\Delta \boldsymbol{f}_k&\leq (\boldsymbol{x}_k-\hat{\boldsymbol{x}}_k)^{\mathrm{T}}\boldsymbol{J}^{\mathrm{T}}\boldsymbol{J}(\boldsymbol{x}_k-\hat{\boldsymbol{x}}_k) \nonumber\\
		\boldsymbol{\Delta f}_k^{\mathrm{T}}\boldsymbol{\Delta f}_k &\leq{\boldsymbol{e}}_{k}^\mathrm{T}\boldsymbol{J}^\mathrm{T}\boldsymbol{J}{\boldsymbol{e}}_k \nonumber \\
		\label{eq32}\beta &\geq 0
	\end{align}
	According to (\ref{eq31}) and (\ref{eq32}), (\ref{eq30}) becomes
	\begin{align}\label{eq33}
		\Delta V_k + \lambda_l V_k - \gamma_l \boldsymbol{w}_k^\mathrm{T}\boldsymbol{w}_k-\gamma_l \boldsymbol{v}_k^\mathrm{T}\boldsymbol{v}_k  \leq 0
	\end{align}
	When $\boldsymbol{w}_k=0$ and $\boldsymbol{v}_k=0$, we obtain $\Delta V_k \leq -\lambda_l V_k < 0$, which proves the global asymptotic stability of the dynamics of the estimation error (\ref{eq19}).
	Moreover, (\ref{eq33}) can be reformulated as:
	\begin{align}
		V_{k+1} &\leq (1-\lambda_l)V_k + \gamma_l \boldsymbol{w}_k^\mathrm{T}\boldsymbol{w}_k+\gamma_l \boldsymbol{v}_k^\mathrm{T}\boldsymbol{v}_k \nonumber \\
		&\leq (1-\lambda_l)V_k + \gamma_l \overline{\boldsymbol{w}}^{\mathrm{T}}\overline{\boldsymbol{w}}+\gamma_l\overline{\boldsymbol{v}}^{\mathrm{T}} \overline{\boldsymbol{v}}
	\end{align}
	which leads to 
	\begin{align}
		V_{k} & \leq (1-\lambda_l)^kV_0 + \gamma_l\sum_{i=0}^{k-1}(1-\lambda_l)^i(\overline{\boldsymbol{w}}^{\mathrm{T}}\overline{\boldsymbol{w}}+\overline{\boldsymbol{v}}^{\mathrm{T}} \overline{\boldsymbol{v}}) \nonumber\\
		\label{eq35}& <  (1-\lambda_l)^kV_0 + \frac{\gamma_l}{\lambda_l}(\overline{\boldsymbol{w}}^{\mathrm{T}}\overline{\boldsymbol{w}}+\overline{\boldsymbol{v}}^{\mathrm{T}} \overline{\boldsymbol{v}})
	\end{align}
	
	Next, we apply Schur complement lemma to (\ref{eq21}), which can be rewritten as
	\begin{align}\label{eq36}
		\boldsymbol{\Upsilon }-\mu^{-1}
		\left[ \begin{matrix}
			\boldsymbol{I}_{n_x} \\	\boldsymbol{0}\\\boldsymbol{0}
		\end{matrix} \right] \left[ \begin{matrix}
			\boldsymbol{I}_{n_x} & 	\boldsymbol{0}& \boldsymbol{0}
		\end{matrix} \right] \succ \boldsymbol{0}
	\end{align} 
	where $\boldsymbol{\Upsilon }= \mathrm{diag}([\lambda_l\boldsymbol{P},(\mu_w-\gamma_l)\boldsymbol{I}_{nw},(\mu_v-\gamma_l)\boldsymbol{I}_{nv} ])$. Pre- and post-multiplying (\ref{eq36}) by $[\boldsymbol{e}_k^\mathrm{T} \quad \overline{\boldsymbol{w}}^\mathrm{T}\quad \overline{\boldsymbol{v}}^\mathrm{T}]$ and its transpose, we have
	\begin{align}\label{eq37}
		\boldsymbol{e}_k^{\mathrm{T}}\boldsymbol{e}_k \leq \mu_l \lambda_l V_k + \mu_l ((\mu_w-\gamma_l)\overline{\boldsymbol{w}}^{\mathrm{T}}\overline{\boldsymbol{w}} + (\mu_v-\gamma_l)\overline{\boldsymbol{v}}^{\mathrm{T}}\overline{\boldsymbol{v}})
	\end{align}
	Finally, substituting (\ref{eq35}) into (\ref{eq37}) yields
	\begin{align}
		\boldsymbol{e}_k^{\mathrm{T}}\boldsymbol{e}_k \leq \mu_l \lambda_l (1-\lambda_l)^kV_0 + \mu_l (\mu_w\overline{\boldsymbol{w}}^{\mathrm{T}}\overline{\boldsymbol{w}} + \mu_v\overline{\boldsymbol{v}}^{\mathrm{T}}\overline{\boldsymbol{v}})
	\end{align}
	which completes the proof.

\begin{remark} 
	Any real matrix $\boldsymbol{J} \in [\underline{\boldsymbol{J}}, \overline{\boldsymbol{J}}] \subseteq \mathbb{R}^{n_x \times n_x}$ can be equivalently represented as a convex combination of vertex matrices:
	\begin{equation}\label{eq39}
		\begin{adjustbox}{max width=\linewidth}
			$
			\boldsymbol{J} = \sum_{i=1}^{N}\boldsymbol{\tau}(i) \boldsymbol{J}_i, \
			\sum_{i=1}^{N}\boldsymbol{\tau}(i)=1, \
			\boldsymbol{\tau}(i) \geq 0, \ N = 2^{n_x^2-\phi}, 
			$
		\end{adjustbox}
	\end{equation}
	where $\boldsymbol{J}_i$ is composed of elements selected from the elements of $\underline{\boldsymbol{J}}$ and $\overline{\boldsymbol{J}}$, i.e. $\boldsymbol{J}_i(j_1,j_2) = \left\lbrace \overline{\boldsymbol{J}}(j_1,j_2), \underline{\boldsymbol{J}}(j_1,j_2) \right\rbrace $, $\phi$ represents the number of corresponding equal elements in  $\underline{\boldsymbol{J}}$ and $\overline{\boldsymbol{J}}$.
	Thus, by substituting (\ref{eq39}) into (\ref{eq20}) and employing the S-procedure, addressing the set-valued problem (related to the presence of the interval variable $\boldsymbol{J}$) in (\ref{eq20}) is equivalent to solving:
	\begin{align}\label{eq40}
		&\quad For \ all \ i = 1,\dots,N  \nonumber \\
		&\left[ \begin{matrix}
			\boldsymbol{\Phi}_{11,i} & * & * & * & * & *\\
			-a\boldsymbol{I} & -b\boldsymbol{I}_{n_x} &* &* & * & *\\
			\boldsymbol{0} & \boldsymbol{0} & -\gamma_l\boldsymbol{I}_{n_w} & * & * & *\\
			\boldsymbol{0} & \boldsymbol{0} & \boldsymbol{0} & -\gamma_l\boldsymbol{I}_{n_v} & * & *\\
			\boldsymbol{\Phi}_{51}& \boldsymbol{P} &\boldsymbol{PD}_1 & -\boldsymbol{WD}_2 & -\boldsymbol{P} & *\\
			\boldsymbol{UJ}_i & \boldsymbol{0} & \boldsymbol{0} & \boldsymbol{0}& \boldsymbol{0} & \boldsymbol{\Phi}_{66}
		\end{matrix}
		\right] \prec \boldsymbol{0}, 
	\end{align}
	where $\boldsymbol{\Phi}_{11,i}=\lambda_l\boldsymbol{P}-\boldsymbol{P}+a(\boldsymbol{J}_i+\boldsymbol{J}_i^\mathrm{T})$. 
	Thus, the optimization problem introduced in (\ref{eq24}) can be transformed to the following one:
	\begin{align}\label{eq41}
		\min_{a, b, \gamma_l, \mu_w, \mu_v, \boldsymbol{P}, \boldsymbol{W}, \boldsymbol{U}} \mu_l \quad \text{s. t.} \quad (\ref{eq40}),(\ref{eq20}).
	\end{align}
	Notice that similar procedures have been applied in former studies \cite{witczak2016lmi, buciakowski2017bounded, } to solve different kind of problems. 
\end{remark}

Based on (\ref{eq23}) , we have
\begin{align}
	\|\boldsymbol{e}_k \|_2 \leq \varphi_k^2 \Rightarrow |\boldsymbol{e}_k(i)| \leq \varphi_k, i = 1,\dots,n_x
\end{align}
where $\varphi_k = \sqrt{\mu_l (\lambda_l(1-\lambda_l)^kV_0 +\mu_w\overline{\boldsymbol{w}}^{\mathrm{T}}\overline{\boldsymbol{w}}+\mu_v\overline{\boldsymbol{v}}^{\mathrm{T}} \overline{\boldsymbol{v}})}$. 
Therefore, a lower bound $\underline{\boldsymbol{x}}_k^p$ and an upper bound $\overline{\boldsymbol{x}}_k^p$ can be computed, based on the peak-to-peak performance index, as follows:
\begin{align}\label{eq43}
	\left\lbrace \begin{aligned}
		&\underline{\boldsymbol{x}}_k^p = \hat{\boldsymbol{x}}_k - \varphi_k \boldsymbol{I}_{n_x}, \\
		&\overline{\boldsymbol{x}}_k^p = \hat{\boldsymbol{x}}_k +  \varphi_k \boldsymbol{I}_{n_x},\ 
	\end{aligned}\right. \ i = 1,\dots,n_x.
\end{align}
Although (\ref{eq43}) can provide an offline interval estimation, 
it is less accurate (tight) due to the conservatism introduced by the overestimation of the effect of the  nonlinear term $\Delta f_k$ done by the inequalities in (31)-(33).
Thus, we apply zonotope analysis next to obtain tighter interval estimation.
\subsection{Online Interval Estimation Based on Zonotope Analysis}
Once the observer in (\ref{eq17}) is established, zonotope analysis is used to compute an interval estimation containing $\boldsymbol{x}_k$  based on the reachable set of  (\ref{eq19}) as stated in the following theorem.
\begin{theorem}
	Given $\boldsymbol{x}_0 \in \mathcal{X}_0 = \langle \boldsymbol{p}_0, \boldsymbol{M}_{\boldsymbol{x},0}\rangle$ and $\hat{\boldsymbol{x}}_0=\boldsymbol{p}_0$, the actual state vector $\boldsymbol{x}_k$ is included inside the zonotope $ \in \mathcal{X}_k = \langle \hat{\boldsymbol{x}}_k, \boldsymbol{M}_{x,k} \rangle$, and its corresponding interval estimation $[\underline{\boldsymbol{x}}_k^z,\overline{\boldsymbol{x}}_k^z]$ can be calculated as follows:
	\begin{align}
		&\boldsymbol{M}_{x,k} = [\tilde{\boldsymbol{A}}\boldsymbol{M}_{x,k-1} \  \boldsymbol{T}_{k-1} \ \boldsymbol{D}_1\boldsymbol{M}_w \ \boldsymbol{LD}_2 \boldsymbol{M}_v ],\\
		&\boldsymbol{T}_{k-1} = [\boldsymbol{D}_{[\boldsymbol{Ja}]_{k-1}}^c \ \boldsymbol{D}_{[\boldsymbol{Ja}]_{k-1}}^s], \\
		& [\boldsymbol{Ja}]_{k-1} = [\boldsymbol{A}]_{k-1}\boldsymbol{M}_{x,k-1},\\
		& [\boldsymbol{A}]_{k-1}=\frac{\partial \boldsymbol{f}}{\partial \boldsymbol{x}}(\mathcal{X}_{k-1}),\\
		\label{eq48}&\left\lbrace \begin{aligned}
			&\underline{\boldsymbol{x}}_k^z(i) = \hat{\boldsymbol{x}}_k - \|\boldsymbol{M}_{x,k}(i,:) \|_1, \\
			&\overline{\boldsymbol{x}}_k^z(i) = \hat{\boldsymbol{x}}_k + \|\boldsymbol{M}_{x,k}(i,:) \|_1 
		\end{aligned}\right. ,\ i = 1,\dots,n_x, \\
		&\boldsymbol{M}_{x,k} = \downarrow_q(\boldsymbol{M}_{x,k})
	\end{align} 
\end{theorem}
where the interval matrix $[\boldsymbol{A}]_{k-1}$ is computed by substituting the zonotopes $\mathcal{X}_{k-1}$ in $\frac{\partial \boldsymbol{f}}{\partial \boldsymbol{x}}(\mathcal{X}_{k-1})$ by their interval outer approximations and then by applying interval arithmetic.
The principles of interval arithmetic are detailed in \cite{jaulin2001interval} and can be implemented using tools such as Cora \cite{althoff2016implementation}.
$\boldsymbol{D}_{[\boldsymbol{Ja}]_{k-1}}^c $ and $ \boldsymbol{D}_{[\boldsymbol{Ja}]_{k-1}}^s$ are matrices related to $[\boldsymbol{Ja}]_{k-1}$ and can be computed by (\ref{eq10}). $\downarrow_q(\cdot)$ is used to limit the column number of $\boldsymbol{M}_{x,k}$.

	First, based on the Lagrange mean value theorem, $\boldsymbol{f}(\boldsymbol{x}_{k-1})$ can be expanded around $\hat{\boldsymbol{x}}_{k-1}$ as
	\begin{align}\label{eq50}
		\boldsymbol{f}(\boldsymbol{x}_{k-1})= \boldsymbol{f}(\hat{\boldsymbol{x}}_{k-1}) + \frac{\partial \boldsymbol{f}}{\partial \boldsymbol{x}}(\boldsymbol{x}'_{k-1})(\boldsymbol{x}_k-\hat{\boldsymbol{x}}_{k-1})
	\end{align}
	where $\boldsymbol{x}'_{k-1}$ is a point between $\boldsymbol{x}_{k-1}$ and $\hat{\boldsymbol{x}}_{k-1}$, and it is contained in $\mathcal{X}_{k-1}$. 
	
	Substituting (\ref{eq50}) into (\ref{eq19}), we have 
	\begin{align}\label{eq51}
		\boldsymbol{e}_{k} = \tilde{\boldsymbol{A}}(\boldsymbol{x}_{k-1}-\hat{\boldsymbol{x}}_{k-1})+ \frac{\partial \boldsymbol{f}}{\partial \boldsymbol{x}}(\boldsymbol{x}'_{k-1})(\boldsymbol{x}_{k-1}-\hat{\boldsymbol{x}}_{k-1}) \nonumber \\+ \boldsymbol{D}_1\boldsymbol{w}_{k-1}-\boldsymbol{L}\boldsymbol{D}_2\boldsymbol{v}_{k-1}
	\end{align}
	Replacing the variables in (\ref{eq50}) by their zonotopic set counterparts leads to,
	\begin{align}\label{eq52}
		\boldsymbol{e}_k \in \tilde{\boldsymbol{A}}(\mathcal{X}_{k-1}-\hat{\boldsymbol{x}}_{k-1}) \oplus \frac{\partial \boldsymbol{f}}{\partial \boldsymbol{x}}(\mathcal{X}_{k-1})(\mathcal{X}_{k-1}-\hat{\boldsymbol{x}}_{k-1}) \nonumber \\
		\oplus\boldsymbol{D}_1\mathcal{W}\oplus\boldsymbol{LD}_2\mathcal{V}
	\end{align} 
	
	Based on Property \ref{Pro1}, (\ref{eq52}) is equivalent to
	\begin{align}
		\boldsymbol{e}_{k} \in& \tilde{\boldsymbol{A}}(\langle \hat{\boldsymbol{x}}_{k-1}, \boldsymbol{M}_{x,k-1}\rangle-\hat{\boldsymbol{x}}_{k-1}) \nonumber \\
		&\oplus \frac{\partial \boldsymbol{f}}{\partial \boldsymbol{x}}(\mathcal{X}_{k-1})(\langle \hat{\boldsymbol{x}}_{k-1}, \boldsymbol{M}_{x,k-1}\rangle-\hat{\boldsymbol{x}}_{k-1}) \nonumber \\
		&\oplus \boldsymbol{D}_1\langle {\boldsymbol{0}}, \boldsymbol{M}_{w}\rangle 
		\oplus \boldsymbol{LD}_2\langle {\boldsymbol{0}}, \boldsymbol{M}_{v}\rangle \nonumber\\
		\in& \tilde{\boldsymbol{A}}\langle {\boldsymbol{0}}, \boldsymbol{M}_{x,k-1}\rangle \oplus [\boldsymbol{Ja}]_{k-1}\langle {\boldsymbol{0}},  \boldsymbol{M}_{x,k-1}\rangle \nonumber \\
		&\oplus \boldsymbol{D}_1\langle {\boldsymbol{0}}, \boldsymbol{M}_{w}\rangle 
		\oplus \boldsymbol{LD}_2\langle {\boldsymbol{0}}, \boldsymbol{M}_{v}\rangle \nonumber\\
		\in& \langle {\boldsymbol{0}}, \tilde{\boldsymbol{A}}\boldsymbol{M}_{x,k-1}\rangle \oplus\langle {\boldsymbol{0}},  [\boldsymbol{Ja}]_{k-1}\boldsymbol{M}_{x,k-1}\rangle \nonumber \\
		\label{eq54}&\oplus \langle {\boldsymbol{0}}, \boldsymbol{D}_1\boldsymbol{M}_{w}\rangle 
		\oplus \langle {\boldsymbol{0}}, \boldsymbol{LD}_2 \boldsymbol{M}_{v}\rangle
	\end{align}
	
	Then, we apply Lemma \ref{lemma2} to $\langle {\boldsymbol{0}},  [\boldsymbol{Ja}]_{k-1}\boldsymbol{M}_{x,k-1}\rangle$, which leads to 
	\begin{align}\label{eq55}
		\lozenge \langle {\boldsymbol{0}},  [\boldsymbol{Ja}]_{k-1}\boldsymbol{M}_{x,k-1}\rangle = \langle {\boldsymbol{0}},  \boldsymbol{T}_{k-1}\rangle
	\end{align}
	where $\boldsymbol{T}_{k-1} = [\boldsymbol{D}_{[\boldsymbol{Ja}]_{k-1}}^c \ \boldsymbol{D}_{[\boldsymbol{Ja}]_{k-1}}^s]$.
	Substituting (\ref{eq55}) into (\ref{eq54}), we obtain
	\begin{align}\label{eq56}
		\boldsymbol{e}_k &\in \langle \boldsymbol{0}, \boldsymbol{M}_{x,k}\rangle \\
		\boldsymbol{M}_{x,k} &= [\tilde{\boldsymbol{A}}\boldsymbol{M}_{x,k-1} \  \boldsymbol{T}_{k-1} \ \boldsymbol{D}_1\boldsymbol{M}_w \ \boldsymbol{LD}_2 \boldsymbol{M}_v ]
	\end{align}
	
	Last, according to (\ref{eq18}) and (\ref{eq56}), it follows
	\begin{align}
		\begin{aligned}
			\boldsymbol{x}_{k} = \hat{\boldsymbol{x}}_k+\boldsymbol{e}_k 
			\in \langle \hat{\boldsymbol{x}}_k, \boldsymbol{M}_{x,k} \rangle \subseteq [\underline{\boldsymbol{x}}_k^z,\overline{\boldsymbol{x}}_k^z]
		\end{aligned}
	\end{align}
	The lower bound $\underline{\boldsymbol{x}}_k^z$ and upper bound $\overline{\boldsymbol{x}}_k^z$ can be computed by (\ref{eq6}), which completes the proof. 

\subsection{Comprehensive method}
Combining the interval estimation (\ref{eq43}) provided by bounded-error and the interval estimation (\ref{eq48}) obtained by zonotope analysis, the bounds of the comprehensive interval estimation are  computed  as follows:
\begin{align}\label{eq57}
	\left\lbrace \begin{aligned}
		&\underline{\boldsymbol{x}}_k = \max(\underline{\boldsymbol{x}}^p_k,\underline{\boldsymbol{x}}^z_k), \\
		&\overline{\boldsymbol{x}}_k = \min(\overline{\boldsymbol{x}}^p_k,\overline{\boldsymbol{x}}^z_k),\ 
	\end{aligned}\right.
\end{align}
where $\underline{\boldsymbol{x}}_k \in \mathbb{R}^{n_x}$ and  $\overline{\boldsymbol{x}}_k\in \mathbb{R}^{n_x}$ are the lower bound  and upper bound of state $\boldsymbol{x}_k$, respectively.
\begin{remark}
	The proposed method also covers linear systems with uncertain matrix $\boldsymbol{A} \in \boldsymbol{\hat{A}} + [\Delta\boldsymbol{A}]$ via $\boldsymbol{f}(\boldsymbol{x}_k) \triangleq \Delta\boldsymbol{A}\boldsymbol{x}_k$, where $\boldsymbol{\hat{A}}$ is a real matrix and $[\Delta\boldsymbol{A}]$ is an interval matrix.
	If $\Delta\boldsymbol{A} \in [\underline{\boldsymbol{J}}, \overline{\boldsymbol{J}}]$, Assumption 1 holds with $\frac{\partial \boldsymbol{f}}{\partial \boldsymbol{x}} = \Delta\boldsymbol{A}$, requiring no extra design.
\end{remark}
\section{Illustrative example}\label{sec4}

The effectiveness of the proposed interval estimation method is demonstrated through the following example:

Consider the nonlinear pendulum model in \cite{tahir2021synthesis}:
\begin{align}
	\left\lbrace \begin{aligned}
		\boldsymbol{x}_{k+1} &= \boldsymbol{Ax}_k+\boldsymbol{f}(\boldsymbol{x}_k)+\boldsymbol{D}_1\boldsymbol{w}_k \\
		\boldsymbol{y}_k &= \boldsymbol{Cx}_k + \boldsymbol{D}_2\boldsymbol{v}_k
	\end{aligned}\right. 
\end{align}
where  
\begin{align}
	\boldsymbol{A}&=\left[\begin{matrix}
		1 & 0.065\\
		0 & 1 
	\end{matrix}\right], \
	\boldsymbol{f}(\boldsymbol{x}_k)=\left[\begin{matrix}
		0 \\
		-0.065\sin(\boldsymbol{x}_k(1)) 
	\end{matrix}\right], \nonumber\\
	\boldsymbol{C}&=[\begin{matrix}
		1 & 0
	\end{matrix}], \ \boldsymbol{D}_1=\boldsymbol{I}_2, \ \boldsymbol{D}_2=0.1,\nonumber 
\end{align}
$\boldsymbol{x}=[x_1 \ x_2]^\mathrm{T}$ denotes the position and angular velocity of pendulum.
The process disturbance and measurement noise satisfy 
\begin{align}
	\boldsymbol{w}_k \in0.006 \left[ \begin{matrix}
		\mathrm{rand}(-1,1)\\
		\mathrm{rand}(-1,1)
	\end{matrix}\right] , \	\boldsymbol{v}_k \in 0.01\mathrm{rand}(-1,1) \nonumber
\end{align}
with $\mathrm{rand}(-1,1)$ means an independent random signal uniformly distributed over $[-1,1]$.
The lower bound $\underline{\boldsymbol{J}}$ and upper bound $\overline{\boldsymbol{J}}$ of the Jacobian matrix are
$
\underline{\boldsymbol{J}}=\left[ \begin{matrix}
	0 & 0\\ -0.065 & 0
\end{matrix}\right],\    
\overline{\boldsymbol{J}}=\left[ \begin{matrix}
	0 & 0\\ 0.065 & 0
\end{matrix}\right] \nonumber
$.
Based on (\ref{eq39}), the number of $\boldsymbol{J}_i$ can be computed as $N=2^{2^2-3}=2$.
Then, setting $\lambda = 0.5$ and solving the optimization problem (\ref{eq41}), we obtain $ \mu_w=15.7, \ \mu_v=15.7, \ \varphi_0 = 1.15$ and  $\boldsymbol{L}=[1.503 \ 7.737]^{\mathrm{T}}$.
Moreover, we set the initial state at $\boldsymbol{x}_0=[1,\ 0]$,  the initial state vector of the observer $\hat{\boldsymbol{x}}_0=[0.98,\ 0.02]$, the shape matrices of the considered zonotopes are $\boldsymbol{M}_{x,0} = \mathrm{diag}([0.1 \ 1]) $, $\boldsymbol{M}_w = 0.006\boldsymbol{I}_2$, $\boldsymbol{M}_v = 0.01$ and  $q = 20$.

\begin{figure}[htbp]
	\begin{center}
		\includegraphics[width=0.82\linewidth]{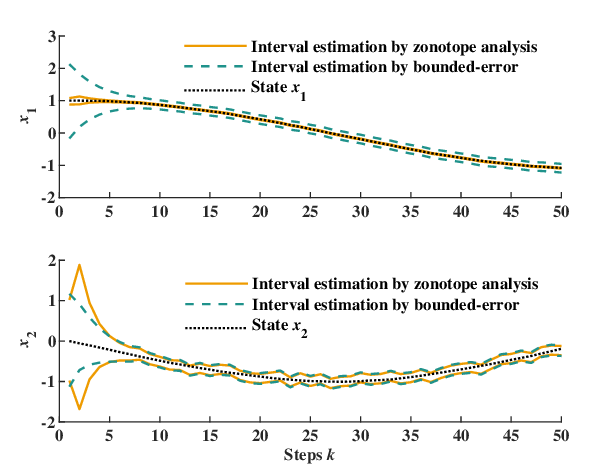}    
		\caption{The interval estimation by zonotope analysis and bounded-error} 
		\label{fig:1}
	\end{center}
\end{figure}
\begin{figure}[htbp]
	\begin{center}
		\includegraphics[width=0.82\linewidth]{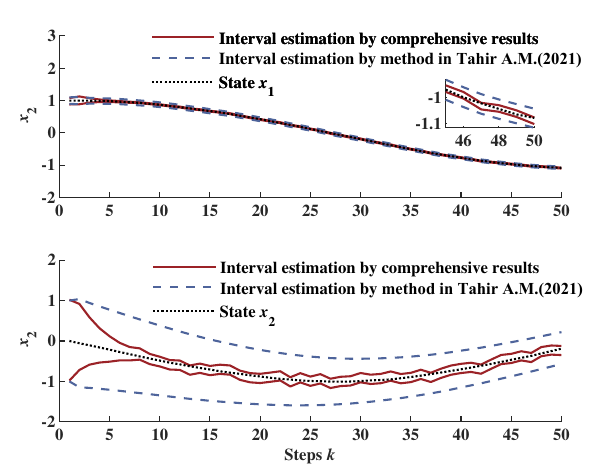}    
		\caption{The interval estimation by comprehensive results and existing method} 
		\label{fig:2}
	\end{center}
\end{figure}

The obtained simulation results by the proposed set-valued state estimators are shown in Fig. \ref{fig:1}. In Fig. \ref{fig:1}, the dotted lines stand for the actual state. The solid lines and the dashed ones correspond to the interval estimation by zonotope analysis (Eq. (\ref{eq48})) and bounded-error method (Eq. (\ref{eq43})), respectively. 
Thanks to the use of inequality relaxations in the observer design process,
zonotope analysis can provide tighter interval estimation than the bounded-error approach.
However, one observes that for the unmeasured state variable, the convergence rate of the bounded-error method is faster than the zonotope analysis approach. 

The comprehensive results (Eq. (\ref{eq57})) make a good balance between estimation accuracy and convergence rate.  The interval estimation obtained by comprehensive approach and the method in \cite{tahir2021synthesis}  is presented in Fig. \ref{fig:2}. In Fig. \ref{fig:2}, the dotted lines stand for the actual state. The solid lines and the dashed ones correspond to the interval estimation by comprehensive approach and the method in \cite{tahir2021synthesis}. 
The simulation results indicate that the proposed method outperforms the existing approach \cite{tahir2021synthesis}. 
In particular, it achieves quicker convergence and provides a more precise estimation interval. 

\section{CONCLUSIONS}\label{sec5}
A comprehensive interval estimation method is proposed for discrete-time bounded Jacobian nonlinear systems subjected to unknown but bounded disturbance and noise. First, a bounded-error method is developed based on a robust observer, which is designed to guarantee the peak-to-peak estimation performance index. Then, a comprehensive interval estimation method is proposed by combining bounded-error information and zonotopic interval estimation obtained from extending zonotope set computation to the dynamics of the estimation error. Finally, simulation results have shown the effectiveness of the developed approach. Extending the proposed method to achieve the fault detection and isolation will be our future works.

\section*{APPENDIX}

\subsection{Size reduction operator}

Let $\overrightarrow{\boldsymbol{M}}$ represent the matrix formed by rearranging the columns of $\boldsymbol{M} \in \mathbb{R}^{n \times m}$ in descending order based on their Euclidean norm. Define $\downarrow_{q}(\boldsymbol{M}) = [ \boldsymbol{M}^1, \boldsymbol{M}^2 ]$, where $\boldsymbol{M}^1$ consists of the first $q-n$ columns of $\overrightarrow{\boldsymbol{M}}$, and $\boldsymbol{M}^2$ is a diagonal matrix with diagonal entries given by 
\begin{align}
	\boldsymbol{M}^2(i, i) = \sum_{j=q-n+1}^{m} \overrightarrow{\boldsymbol{M}}(i, j), \ i = 1, \dots, n. \nonumber
\end{align}
Note that the reduce operation limits the column number of the generator matrix, but introduces additional conservatism.
\bibliography{ifacconf}             
                                                   







\end{document}